\documentclass[12pt,a4paper]{article}
\usepackage[french]{babel}
\usepackage{amsfonts}
\newtheorem{exe}{Exercice}
\newtheorem{Theo}{Th\'eor\`eme}[section]
\newtheorem{lemm}{Lemme}[section]

\newcommand{\dint}{\displaystyle\int}
\newcommand{\dfrac}{\displaystyle\frac}

\newcommand\dlimsup{\displaystyle\limsup}

\newcommand{\é}{\'e}

\newcommand{\à}{\`a}

\begin{document}
\title{\bf {Identification d'un processus autor\égressif gaussien
stable par la m\éthode de moyennisation logarithmique  dans le cas
r\éel}}
\author{\bf {Faouzi Chaabane}\thanks{Equipe d'Analyse Stochastique et Mod\élisation Statistique, Facult\é des Sciences de Bizerte, 7021 Jarzouna (Tunisie). Tel:+21672590613; Fax:+21672590566, E-mail address:faouzi.chaabane@fsb.rnu.tn.} et \bf {Hamdi Fathallah} \thanks{Laboratoire LMV, Universit\é de Versailles Saint-Quentin-En-Yvelines, 45 Avenue des Etats-Unis Batiment Ferma 78035 Versailles (France). Tel:+33139253629; Fax:+33139254645, E-mail address:hamdi.fathallah@math.uvsq.fr}}
\date{juin 2006}
\maketitle
\begin{abstract}
Dans ce travail, on consid\ère un mod\èle autor\égressif gaussien
stable \à temps continu unidimentionnel et on lui applique les
th\éor\èmes limites par moyennisation logarithmique obtenus pour des
martingales locales continues \à temps continu. On construit alors
un estimateur de la covariance du bruit $\sigma^{2}$ et un autre
estimateur de $\theta$ autre que celui des moindres carr\és. En
exploitant la m\éthode de pond\ération, on am\éliore les vitesses de
convergence de ces nouveaux estimateurs.\\
\begin{center}
{\bf Identification of a stable gaussian autoregression process by
logarithmic averaging method in the real case }
\end{center}
\begin{center}
{\bf Abstract}
\end{center}
\par In the present work, we consider a  stable one-dimensional gaussian
autoregressive model  in continous time. Using the limit theorems
with logarithmic averaging obtained for continous local martingales,
we construct then an estimator of the noise covariance $\sigma^{2}$
and an estimator of $\theta$ different of the one of the least
squares estimator. By exploiting the weighting method
we ameliorate the convergence rates of these new estimators.
 \\Logarithmic central limit theorem; law of iterated logarithm.
\end{abstract}
\section{Introduction}
Le but de ce travail est d'estimer les param\ètres d'un mod\èle
autor\égressif gaussien stable \à temps continu unidimentionnel.
Etant donn\é un mouvement brownien standard r\'eel $B=(B_{t},\,
t\geq 0)$. On d\'efinit le processus autor\égressif
$X=(X_{t},t\geq0)$ avec $X_{0}=0$ par la relation:
\begin{equation}\label{e1}
X_{t}= \theta \dint_{0}^{t}X_{s}ds+ \sigma B_{t};\quad t\geq 0,
\end{equation}
avec $\theta$ $\leq 0$ \,et\, $\sigma$ un param\`etre r\'eel.\\On
d\ésigne par $ \tilde{\theta}_{t}$  l'estimateur des moindres
carr\'es pond\ér\é  de $ \theta$ d\'efini par
\begin{equation}\label{e2}
\tilde{\theta}_{t}=P_{t}^{-1}\dint_{0}^{t}\omega_{s}X_{s}dX_{s},\quad
t\geq 0,
\end{equation}
correspondant au poids $(\omega_{s})$ donn\é par
\begin{equation}\label{e8}
\omega_{s}=s^{-\frac{\alpha}{2}}\exp\{\frac{s^{1-\alpha}}{2(1-\alpha)}\};\quad
\dfrac{1}{2} < \alpha < 1,
\end{equation}
 avec $P_{t}=\dint_{0}^{t}\omega_{s}X_{s}^{2}ds$
et par
$\bar{\theta}_{t}=\frac{1}{t}\dint_{0}^{t}\tilde{\theta}_{s}ds$\,
son
moyennis\é.\\
   On remarque que pour $\omega_{s}=1$  on obtient l'estimateur des moindres carr\'es
 $\widehat{\theta}_{t}$ de $ \theta$  d\'efini par
\begin{equation}\label{e3}
\widehat{\theta}_{t}= ( \zeta_{t})^{-1}\int_{0}^{t}X_{s}
dX_{s},\quad t\geq 0,\quad \mbox{avec}\quad
\zeta_{t}=\int_{0}^{t}X_{s}^{2}ds.
\end{equation}
Dans \cite{DA}, dans le cas multidimentionnel, Darwich montre que
cet estimateur v\érifie la loi du logarithme it\ér\é (LLI). De
fa\c{c}on pr\'ecise, on a
$$ |\widehat{\theta}_{t}- \theta |= \mathcal{O}\bigg( (\frac{\log\log t}{t})^{\frac{1}{2}}\bigg)\quad p.s..$$
Ce r\ésultat, dans le cas unidimentionnel, d\écoule de la LLI pour
les martingales r\éelles continues (voir \cite{PL}). Les th\éor\èmes
limites par moyennisation logarithmique permettent entre autre de
montrer des r\ésultats aussi bien pour l'estimateur
des moindres carr\'es que l'estimateur pond\ér\é de type\\
$\ast$ Th\'eor\`eme de la limite centrale presque-s\^ure\\\\
TLCPS: \quad
 $ \dfrac{1}{\log t}\dint_{0}^{t}\dfrac{ds}{s}\,
\delta_{\{\sqrt{s}( \widehat{\theta}_{s}- \theta)\}} \Longrightarrow
 \mathcal{N}(0,2 \theta)\quad p.s..$\\\\
 Avec ($\Longrightarrow$) d\énote la convergence en loi.\\
La convergence en moyenne d'ordre deux associ\ée au TLCPS permet
d'une part de construire un estimateur de $\sigma$ et d'autre part
de donner un
autre estimateur de $\theta$ \\
$\ast$\quad $\hat{\sigma}_{t}^{2}=\dfrac{1}{\log t}\dint_{0}^{t} (
\widehat{\theta}_{s}-
 \bar{\theta}_{t})^{2}X_{s}^{2}ds \longrightarrow  \sigma^{2}\quad
 ~p.s.,~~(t\longrightarrow\infty).$\\
$\ast$\quad$ \check{\theta}_{t}=\dfrac{1}{2\log t}\dint_{0}^{t}
(\widehat{\theta}_{s}- \bar{\theta}_{t})^{2}ds \longrightarrow
\theta\quad
 ~p.s.,~~(t\longrightarrow\infty).$\\
 On donnera par la suite les vitesses de convergence en loi et presque s\^ure
de l'estimateur $\check{\theta}_{t}$ de $\theta$. La m\éthode de
pond\ération nous permettera d'am\éliorer les vitesses de
convergence de ces estimateurs. les principaux r\ésultats seront
\énonc\és au paragraphe $2$ et leurs preuves seront donn\ées au
paragraphe $4$. Le paragraphe $3$ sera consacr\é \à la donn\ée de
quelques outils de d\émonstration.
\section{Enonc\'e des principaux r\'esultats}
\subsection{Cas sans pond\ération $(\omega_{s}=1)$}
\begin{Theo}\label{t2}
Soit $X=(X_{t},\, t\geq 0)$ le r\égresseur stochastique d\éfini par
{\normalfont{(\ref{e1})}}. Alors on a les r\'esultats suivants
\begin{enumerate}
\item Th\'eor\`eme de la limite centrale presque-s\^ure\\\\
TLCPS: \quad $ \dfrac{1}{\log t}\dint_{0}^{t}\dfrac{ds}{s}\,
\delta_{\{\sqrt{s}( \widehat{\theta}_{s}- \theta)\}} \Longrightarrow
 \mathcal{N}(0,2 \theta)\quad p.s..$\\
\item La loi forte quadratique\\\\
   i)\quad $ \dfrac{1}{\log t}\dint_{0}^{t} \zeta_{s}^{2}( \widehat{\theta}_{s}-
 \theta)^{2}\dfrac{ds}{s^{2}}\,
 \longrightarrow \dfrac{\sigma^{4}}{2 \theta}\quad ~p.s.,~~
 (t\longrightarrow\infty)$.\\\\
ii) \quad $\hat{\sigma}_{t}^{2}=\dfrac{1}{\log t}\dint_{0}^{t} (
\widehat{\theta}_{s}-
 \bar{\theta}_{t})^{2}X_{s}^{2}ds \longrightarrow  \sigma^{2}\quad
 ~p.s.,~~(t\longrightarrow\infty)$.\\\\
iii)\quad  Si de plus on suppose que les observations $(X_{t},\,
t\geq 0)$ v\'erifient l'hypoth\èse suivante
\begin{equation}\label{e6}
(H1)\quad
\dfrac{1}{t}\int_{0}^{t}X_{s}^{2}ds=\frac{\sigma^{2}}{2\theta}+
\mathbf{o}\bigg((\log\log t)^{-1}\bigg)\quad p.s.,
\end{equation}
on obtient
$$ \check{\theta}_{t}=\frac{1}{2\log t}\int_{0}^{t} (\widehat{\theta}_{s}- \bar{\theta}_{t})^{2}ds \longrightarrow  \theta\quad
 ~p.s.,~~(t\longrightarrow\infty).$$
 \end{enumerate}
 \end{Theo}
 \begin{Theo}\label{t3}
On se place dans le cadre du th\'eor\`eme {\normalfont{\ref{t2}}} et
en renfor\c{c}ant l'hypoth\`ese $(H1)$ de la mani\`ere suivante
\begin{equation}\label{e7}
(H2)\quad
\dfrac{1}{t}\int_{0}^{t}X_{s}^{2}ds=\frac{\sigma^{2}}{2\theta}+
\mathbf{o}\bigg((\log\log t)^{-1}(\log t)^{-\frac{1}{2}}\bigg)\quad
p.s.,
\end{equation}
 on obtient
\begin{enumerate}
\item Th\'eor\`eme de la limite centrale logarithmique\\\\
TLCL:\quad $ \sqrt{\log t}\,\Bigg(\dfrac{1}{2\log t}\dint_{0}^{t}
(\widehat{\theta}_{s}- \bar{\theta}_{t})^{2}ds- \theta \Bigg)
\Longrightarrow \mathcal{N}(0,(2\theta)^{2}).$\\
\item La loi du logarithme it\'er\'e logarithmique\\\\
LLIL: \quad $ \dlimsup_{t\to \infty}\dfrac{ \log
t}{\sqrt{\log\log\log t}}\,\Bigg |\dfrac{1}{2\log t}\dint_{0}^{t}
(\widehat{\theta}_{s}- \bar{\theta}_{t})^{2}ds-\theta \Bigg|=
2\theta \sqrt{2}\quad p.s..$
\end{enumerate}
\end{Theo}
\subsection{Cas avec pond\ération}
\begin{Theo}\label{t4}
Soit $X=(X_{t},\, t\geq 0)$ le processus autor\égressif gaussien \`a
temps continu d\éfini par la relation {\normalfont{(\ref{e1})}}.
Alors l'estimateur de moindre carr\'e pond\ér\é $
\tilde{\theta}_{t}$ de $ \theta$ donn\é par la relation
{\normalfont{(\ref{e2})}} ainsi son moyennis\é $\bar{\theta}_{t}$
convergent au sens presque-s\^ur. De fa\c{c}on pr\'ecise, pour $
\dfrac{1}{2}<\alpha<1$, on a
$$|\tilde{\theta}_{t}- \theta |= \mathcal{O}\bigg( (\frac{\log
t}{t^{\alpha}})^{\frac{1}{2}}\bigg)\quad p.s..$$ Si de plus on
suppose que $2\alpha-1\leq \alpha'< \frac{3}{2}\alpha-\frac{1}{2}$,
on a $$~~|\bar{\theta}_{t}- \theta |= \mathcal{O}\bigg(
(\frac{\log\log t}{t})^{\frac{1}{2}}\bigg)\quad p.s..$$
\end{Theo}
\begin{Theo}\label{t5}
Sous les hypoth\`eses du th\'eor\`eme {\normalfont{\ref{t4}}}, on a
les r\'esultats suivants
\begin{enumerate}
\item Th\'eor\`eme de la limite centrale presque-s\^ure\\\\
TLCPS: \quad $ \dfrac{1-\alpha}{
t^{1-\alpha}}\dint_{0}^{t}\dfrac{ds}{s^{\alpha}}\,
\delta_{\{s^{\frac{\alpha}{2}}( \tilde{\theta}_{s}- \theta)\}}
\Longrightarrow
 \mathcal{N}(0,2 \theta)\quad p.s..$\\
\item La loi forte quadratique\\\\
   i)\quad $\dfrac{1-\alpha}{ t^{1-\alpha}}\dint_{0}^{t} \frac{P_{s}^{2}}{U_{s}^{2}}( \tilde{\theta}_{s}-
 \theta)^{2}ds \,
 \longrightarrow \dfrac{\sigma^{4}}{2 \theta}\quad ~p.s.,~~
 (t\longrightarrow\infty).\\\\$
ii) \quad $\tilde{\sigma}_{t}^{2}=\dfrac{4(1-\alpha)}{
t^{1-\alpha}}\dint_{0}^{t} ( \tilde{\theta}_{s}-
 \bar{\theta}_{t})^{2}X_{s}^{2}ds \longrightarrow  \sigma^{2}\quad
 ~p.s.,~~(t\longrightarrow\infty).$\\\\
iii)\quad  Si de plus on suppose que les observations  $(X_{t},
\,t\geq 0)$ v\érifient l'hypoth\èse  suivante   $$(H3)\quad
\quad\quad
\frac{1}{t}\dint_{0}^{t}X_{s}^{2}ds-\frac{\sigma^{2}}{2\theta}=o(t^{\alpha'-1})
~p.s.,~~ (t\longrightarrow\infty),\quad \mbox{avec}\quad
\frac{1}{2}\leq\alpha'<\alpha,$$
 on d\égage un estimateur fortement
consistant de $\theta$ \à savoir
$$ \breve{\theta}_{t}=\frac{1-\alpha}{ 2t^{1-\alpha}}\int_{0}^{t} (\tilde{\theta}_{s}- \bar{\theta}_{t})^{2}ds \longrightarrow  \theta\quad
 ~p.s.,~~(t\longrightarrow\infty).$$
\end{enumerate}
\end{Theo}
\begin{Theo}\label{t6}
Soit $X=(X_{t},\, t\geq 0)$ le processus autor\égressif gaussien \à
temps continu satisfisant l'\équation {\normalfont{(\ref{e1})}}.
Supposons que les observations  $(X_{t}, \,t\geq 0)$ v\érifient
l'hypoth\èse suivante
$$(H4)\quad \quad\quad \frac{1}{t}\dint_{0}^{t}X_{s}^{2}ds-\frac{\sigma^{2}}{2\theta}=o(t^{\alpha'-1})
~p.s.,~~ (t\longrightarrow\infty),\quad \mbox{pour}\quad
\frac{1}{2}\leq \alpha'< 3\alpha-2.$$ Alors, on obtient
\begin{enumerate}
\item Th\'eor\`eme de la limite centrale logarithmique\\\\
TLCL:\quad $  t^{\frac{1-\alpha}{2}}\,\Bigg (\dfrac{1-\alpha}{
2t^{1-\alpha}}\dint_{0}^{t} (\tilde{\theta}_{s}-
\bar{\theta}_{t})^{2}ds-\theta \Bigg) \Longrightarrow
\mathcal{N}(0,4\theta^{2}(1-\alpha)).$
\item La loi du logarithme it\'er\'e logarithmique\\\\
LLIL: \quad $ \dlimsup_{t\to \infty} {\dfrac{
t^{1-\alpha}}{{\sqrt{\log\log t^{1-\alpha}}}}}\,\Bigg
|\dfrac{1-\alpha}{ 2t^{1-\alpha}}\dint_{0}^{t} (\tilde{\theta}_{s}-
\bar{\theta}_{t})^{2}ds- \theta \Bigg|=
\sqrt{2(1-\alpha)}2\theta\,\quad p.s..$
\end{enumerate}
\end{Theo}
\section{Les outils des d\'emonstrations}
Au d\ébut de ce paragraphe, on introduit
$\tilde{M}=(\tilde{M}_{t},t\geq 0) $ la martingale locale r\'eelle
continue d\'efinie par $$
\tilde{M}_{t}=\dint_{0}^{t}\omega_{s}X_{s}dB_{s},\quad t\geq 0,$$
dont le processus croissant pr\évisible  $ \langle
\tilde{M}\rangle=( \langle \tilde{M}\rangle_{t}, t\geq 0)$ est
donn\é par
\begin{equation}\label{e4}
\langle\tilde{M}\rangle_{t}=\int_{0}^{t}\omega_{s}^{2}X_{s}^{2}ds.
\end{equation}
D'une part, d'apr\ès (\ref{e1}) et (\ref{e2}), on a
\begin{equation}\label{e16}
\tilde{M}_{t}=\sigma^{-1}P_{t}(\tilde{\theta}_{t}-\theta);\quad
 t\geq 0.
\end{equation}
D'autre part, on a d'apr\`es les relations (\ref{e1}) et (\ref{e3})
 \begin{equation}\label{e5}
 M_{t}=\sigma^{-1}\zeta_{t}(\widehat{\theta}_{t}-\theta);\quad
 t\geq 0,
\end{equation}o\`u $M=(M_{t},t\geq 0)$ est la martingale locale
r\'eelle continue d\'efinie par\\ $$
M_{t}=\int_{0}^{t}X_{s}dB_{s},\quad  t\geq 0,$$ dont le processus
croissant pr\évisible $ \langle M\rangle=(\langle M\rangle_{t},
t\geq 0)$
n'est autre que le processus $ \zeta=( \zeta_{t}, t\geq 0) $ introduit dans la relation (\ref{e3}).\\
Afin de simplifier les preuves des principaux r\'esultats, on
\étudiera les comportements asymptotiques des processus $(\langle M
\rangle_{t}, \,t\geq 0)$, $(\langle \tilde{M}\rangle_{t}, \,t\geq
0)$ et $(P_{t},\, t\geq 0)$. On donnera ensuite quelque propri\ét\és
de la pond\ération $(\omega_{t})$.\\ Le lemme suivant (voir
\cite{c3}) donne Le comportement de la variation quadratique
pr\évisible de la martingale $M$.
\begin{lemm}\label{c}
Soit le processus $(X_{t}, \,t\geq 0)$ d\éfini par (\ref{e1}). Alors
 \begin{equation}\label{e9} \frac{1}{t}\int_{0}^{t}X_{s}^{2}ds\longrightarrow \frac{\sigma^{2}}{2 \theta}\quad
 ~p.s.,~~(t\longrightarrow\infty).
\end{equation}\\
\end{lemm}
Pour  $(\omega_{t})$, le poids d\éfini par
{\normalfont{(\ref{e8})}}, on introduit les deux processus  suivants
$$V_{t}= \Big( \int_{0}^{t}\omega_{s}^{2}ds \Big)^{\frac{1}{2}}\quad
et\quad U_{t}=\int_{0}^{t}\omega_{s}ds,\quad t\geq 0.$$ On a alors
les lemmes suivants
\begin{lemm}\label{a}
Le poids $(\omega_{t})$ satisfait les propri\ét\ées suivantes:\\\\
 $ P_{1}) ~~ t^{-\alpha}\omega_{t}^{-1}U_{t}=2+
 \mathbf{o}(t^{\alpha-1}),~~
(t\longrightarrow\infty).$\\\\
 $ P_{2}) ~~ t^{-\alpha}\omega_{t}^{-2}V_{t}^{2} =1
                  + \mathbf{o}(t^{\alpha-1}),~~
(t\longrightarrow\infty).$ \\\\
 $ P_{3})~~\dfrac{t^{1-\alpha}}{1-\alpha}-\log
V_{t}^{2}=\mathbf{o}(\log{t}),~~
(t\longrightarrow\infty).$\\\\
 $ P_{4})~~\dint_{0}^{t}\frac{V_{s}^{2}}{U_{s}^{2}}ds-\frac{1}{1-\alpha}t^{1-\alpha}=\mathbf{o}(\log{t}),~~
(t\longrightarrow\infty).$
\end{lemm}
\begin{lemm}\label{b}
 On suppose que les observations  $(X_{s}, \,s\geq 0)$ v\'erifient l'hypoth\èse $(H3)$
 du th\éor\ème {\normalfont{\ref{t5}}} suivante
$$\quad \quad\quad \frac{1}{t}\dint_{0}^{t}X_{s}^{2}ds-\frac{\sigma^{2}}{2\theta}=o(t^{\alpha'-1})
~p.s.,~~ (t\longrightarrow\infty)\quad  \mbox{pour}
\,\,\,\frac{1}{2}\leq \alpha'< \alpha.$$
 Alors, on a\\\\
$ i)~~ \dfrac{\langle
\tilde{M}\rangle_{t}}{V_{t}}-\frac{\sigma^{2}}{2\theta}=o(t^{\alpha'-\alpha})~p.s.,~~
(t\longrightarrow\infty).$\\\\
 $ii)~~
\dfrac{P_{t}}{U_{t}}-\frac{\sigma^{2}}{2\theta}=o(t^{\alpha'-\alpha})
~p.s.,~~ (t\longrightarrow\infty).$
\end{lemm}
La preuve de ces deux lemmes est donn\ée dans l'annexe.
\section{D\émonstration des principaux r\'esultats}
\subsection*{Preuve du th\'eor\`eme \ref{t2}}
\begin{enumerate}
\item En appliquant le th\'eor\`eme de la limite centrale presque-s\^ure pour le couple $(M,V)$ avec $V_{t}^{2}=t$ (voir th\éor\ème $1$ dans \cite{Fa}), on obtient$$ (\log
t)^{-1}\int_{0}^{t}\frac{ds}{s}\, \delta_{\{ \frac{M_{s}}{
\sqrt{s}}\}} \Longrightarrow
\mathcal{N}(0,\frac{\sigma^{2}}{2\theta}) \quad p.s..$$ Vu les deux
relations (\ref{e5}) et (\ref{e9}) on a
\begin{equation}\label{e12}
\frac{M_{s}}{ \sqrt{s}}\,\sim \, \frac{ \sigma}{2 \theta}
\sqrt{s}(\widehat{\theta}_{s}- \theta)\quad
p.s.,~~(s\longrightarrow\infty).
\end{equation}
Soit $$ \Delta_{t}= \dint_{0}^{t} \varphi(\frac{M_{s}}{
\sqrt{s}})\frac{ds}{s}- \dint_{0}^{t} \varphi(\frac{ \sigma}{2
\theta} \sqrt{s}(\widehat{\theta}_{s}- \theta))\frac{ds}{s},$$ o\`u
$ \varphi$ est une fonction lipschitzienne  continue.\\ Gr\^ace \`a
l'\'equivalence (\ref{e12}), on a $$(\log t)^{-1}|
\Delta_{t}|\,\longrightarrow \,0\quad
p.s.,~~(t\longrightarrow\infty).$$ Par cons\équent  $$ (\log
t)^{-1}\int_{0}^{t}\frac{ds}{s}\,\delta_{\{ \frac{\sigma}{2\theta}
\sqrt{s}(\widehat{\theta}_{s}- \theta) \}} \Longrightarrow
\mathcal{N}(0,\frac{\sigma^{2}}{2\theta} ).$$ Le r\'esultat en
d\écoule.\\
\item
i)\quad D'apr\`es la relation (\ref{e5}), on a
\begin{equation}\label{e40}
 M_{s}^{2}=\sigma^{-2} \langle M\rangle_{s}^{2} (
\widehat{\theta}_{s}- \theta)^{2}.
\end{equation}
 La loi forte quadratique une (LFQ1) appliqu\'ee \`a la martingale $M$ normalis\'ee par le
processus $(V_{t}= \sqrt{ t},t\geq 0)$ (voir th\éor\ème $3$ dans
\cite{Fa}) donne
\begin{equation}
(\log t)^{-1}\int_{0}^{t}\frac{M_{s}^{2}}{s}\frac{ds}{s}
\longrightarrow \frac{\sigma^{2}}{2 \theta}\quad p.s..
\end{equation}
Compte tenu de la relation (\ref{e40}), on voit que
\begin{equation}\label{e13}(\log
t)^{-1}\int_{0}^{t}\langle M\rangle_{s}^{2}(\widehat{\theta}_{s}-
\theta )^{2}\frac{ds }{s^{2}} \longrightarrow \frac{\sigma^{4}}{2
\theta}\quad p.s..
\end{equation}
 Ainsi le r\ésultat est \établi.\\\\
 ii)\quad En appliquant la loi forte quadratique deux (LFQ2) au couple $(M,V)$ (voir th\éor\ème $3$ dans
\cite{Fa}), on obtient
 $$(\log t)^{-1}\dint_{0}^{t}\frac{M_{s}^{2}}{ \langle M \rangle_{s}^{2} }d\langle M \rangle_{s}\longrightarrow 1\quad
 p.s.,~~(t\longrightarrow\infty). $$
 Gr\^ace \à  la relation (\ref{e5}), on a
  \begin{equation}\label{e25}
 (\log t)^{-1}\dint_{0}^{t} ( \widehat{\theta}_{s}-
 \theta)^{2}X_{s}^{2}ds \,\longrightarrow\, \sigma^{2}\quad
 p.s.,~~(t\longrightarrow\infty).
\end{equation}
D'une part, dans  \cite{DA}, Darwich  a montr\é que
\begin{equation}\label{e39}
|\widehat{\theta}_{t}- \theta| =\mathcal{O}\Bigg((\dfrac{\log\log
t}{t})^{\frac{1}{2}}\Bigg)\quad p.s..
\end{equation}
 D'autre part, vu que
\begin{equation}\label{e27}
\bar{\theta}_{t}-\theta=\dfrac{1}{t}\dint_{0}^{t}(\tilde{\theta}_{s}-\theta)
\,ds,
\end{equation} et en tenant compte de  la relation (\ref{e39}), on obtient
\begin{equation}\label{e24} |\bar{\theta}_{t}- \theta|
=\mathcal{O}\Bigg((\dfrac{\log\log t}{t})^{\frac{1}{2}} \Bigg)\quad
p.s..
\end{equation}
On remarque que
\begin{eqnarray*}
   \bigg|(\log t)^{-1}\dint_{0}^{t}(\widehat{\theta}_{s}- \theta
)^{2}X^{2}_{s}ds-(\log t)^{-1}\dint_{0}^{t}(\widehat{\theta}_{s}-
\bar{\theta}_{t} )^{2}X^{2}_{s}ds\bigg|\leq~~~~~~~~~~~~~~~~~~~~~~~~~~\\\\
\underbrace{(\log t)^{-2}(\bar{\theta}_{t}- \theta
)\dint_{0}^{t}(\widehat{\theta}_{s}- \bar{\theta}_{t}
)X^{2}_{s}ds}_{(C_{t})} +\underbrace{(\log t)^{-1}(\bar{\theta}_{t}-
\theta )^{2}\dint_{0}^{t}X^{2}_{s}ds}_{(D_{t})}.
\end{eqnarray*}
En utilisant les relations (\ref{e9}) et (\ref{e24}), on obtient
$$C_{t}\longrightarrow0\quad
 p.s.\quad \mbox{et} \quad D_{t}\longrightarrow0\quad
 p.s.,~~(t\longrightarrow\infty). $$
Ce qui implique que
\begin{equation}\label{e35}
\bigg|(\log t)^{-1}\dint_{0}^{t}(\widehat{\theta}_{s}- \theta
)^{2}X^{2}_{s}ds-(\log t)^{-1}\dint_{0}^{t}(\widehat{\theta}_{s}-
\bar{\theta}_{t})^{2}X^{2}_{s}ds\bigg|\longrightarrow0\quad
 p.s.,~~(t\longrightarrow\infty).
\end{equation}
Le r\ésultat d\écoule des convergences (\ref{e25}) et
(\ref{e35}).\\\\
iii)\quad La propri\'et\'e (\ref{e13}) s'\'ecrit $$ (\log
t)^{-1}\Bigg( \int_{0}^{t}(\widehat{\theta}_{s}- \theta
)^{2}\bigg(\frac{ \langle M\rangle_{s}^{2}}{s^{2}} -\frac{
\sigma^{4} }{4 \theta^{2}}\bigg)ds+\frac{ \sigma^{4} }{4
\theta^{2}}\int_{0}^{t}(\widehat{\theta}_{s}- \theta )^{2}ds
\Bigg)\,\longrightarrow\,\frac{\sigma^{4}}{2 \theta}\quad
p.s.,~~(t\longrightarrow\infty).$$
 Gr\^ace \à la relation
(\ref{e39}), on obtient
$$ (\log t)^{-1}\int_{t_{0}}^{t}(\widehat{\theta}_{s}- \theta )^{2}\bigg(\frac{
\langle M\rangle_{s}^{2}}{s^{2}} -\frac{ \sigma^{4} }{4
\theta^{2}}\bigg)ds= (\log t)^{-1}\int_{t_{0}}^{t} \mathcal{O}\bigg(
\frac{\log\log s}{s} \bigg)\bigg(\frac{ \langle
M\rangle_{s}^{2}}{s^{2}} -\frac{ \sigma^{4} }{4 \theta^{2}}\bigg)ds.
$$
Vu l'hypoth\`ese (H1), on a
$$(\log t)^{-1}\dint_{0}^{t}(\widehat{\theta}_{s}- \theta
)^{2}\bigg(\frac{ \langle M\rangle_{s}^{2}}{s^{2}}-\frac{ \sigma^{4}
}{4 \theta^{2}}\bigg)ds\, \longrightarrow\,0\quad
p.s.,~~(t\longrightarrow\infty).$$ Par cons\équent
\begin{equation}\label{e33}
  (\log
t)^{-1}\int_{0}^{t}(\widehat{\theta}_{s}- \theta )^{2}ds\,
\longrightarrow\,2 \theta\quad p.s.,~~(t\longrightarrow\infty).
\end{equation}
  Par ailleurs,
on voit que
\begin{eqnarray*}
 \Bigg|(\log
t)^{-1}\dint_{0}^{t}(\widehat{\theta}_{s}- \theta )^{2}ds-(\log
t)^{-1}\dint_{0}^{t}(\widehat{\theta}_{s}- \bar{\theta}_{t}
)^{2}ds\Bigg|\leq~~~~~~~~~~~~~~~~~~~~~~~~~~~~~~~~\\\\\underbrace{t(\log
t)^{-1}(\bar{\theta}_{t}- \theta )^{2}}_{(G_{t})}+\underbrace{2(\log
t)^{-1}(\theta-\bar{\theta}_{t})\int_{0}^{t}(\widehat{\theta}_{s}-
\bar{\theta}_{t} )ds}_{(H_{t})}.\end{eqnarray*}
 D'apr\ès la relation $(\ref{e24})$, on obtient
$$G_{t}\longrightarrow 0 \quad p.s.\quad \mbox{et}\quad H_{t}\longrightarrow 0 \quad p.s.,~~(t\longrightarrow\infty). $$
 Par cons\équent
\begin{equation}\label{e34}
\Bigg|(\log t)^{-1}\dint_{0}^{t}(\widehat{\theta}_{s}- \theta
)^{2}ds-(\log t)^{-1}\dint_{0}^{t}(\widehat{\theta}_{s}-
\bar{\theta}_{t} )^{2}ds\Bigg|\longrightarrow 0\quad
p.s.,~~(t\longrightarrow\infty).
\end{equation}
 Alors, d'apr\ès les
convergences   $(\ref{e33})$ et $(\ref{e34})$, on conclut le
r\ésultat.
 \end{enumerate}
\subsection*{Preuve du th\'eor\`eme \ref{t3}}
\begin{enumerate}
\item Gr\^ace \à l'hypoth\`ese (H2), on a $$\dfrac{\langle M\rangle_{t}}{t}=\frac{\sigma^{2}}{2\theta}+
\mathbf{o}\bigg((\log t)^{-2}\bigg)\quad p.s..$$ Alors en appliquant
le th\éor\ème de la limite centrale logarithmique au couple $(M,V)$
(voir th\éor\ème $5$ dans \cite{Fa}) pour $ V_{t}^{2}=t$ et
$f(x)=x^{2}-1$, on obtient  $$ (\log
t)^{-\frac{1}{2}}\int_{1}^{t}\bigg(
\frac{M_{s}^{2}}{s}-\frac{\sigma^{2}}{2\theta}\bigg)\frac{ds}{s}
 \Longrightarrow \mathcal{N}(0,\frac{\sigma^{4}}{
 \theta^{2}}).$$
 Par suite $$ (\log t)^{-\frac{1}{2}}\int_{1}^{t}
\frac{M_{s}^{2}}{s^{2}}ds- \frac{\sigma^{2}}{2\theta}(\log
t)^{\frac{1}{2}} \Longrightarrow \mathcal{N}(0,\frac{\sigma^{4}}{
 \theta^{2}}).
$$Ainsi vu la relation (\ref{e5}), il vient que
$$(\log t)^{-\frac{1}{2}}\int_{1}^{t}
\frac{M_{s}^{2}}{s^{2}}ds=\sigma^{-2}(\log
t)^{-\frac{1}{2}}\int_{1}^{t}\frac{ \langle
M\rangle_{s}^{2}}{s^{2}}( \widehat{\theta}_{s}- \theta )^{2}ds.$$
D\ésormais, on pose
$$ I_{t}=\int_{1}^{t}\frac{ \langle M\rangle_{s}^{2}
}{s^{2}}(\widehat{\theta}_{s}- \theta)^{2} ds.$$Donc
\begin{equation}\label{e14}
(\log t)^{-\frac{1}{2}}I_{t}=\int_{1}^{t}(\widehat{\theta}_{s}-
\theta)^{2}\bigg( \frac{ \langle M\rangle_{s}^{2}
}{s^{2}}-\frac{\sigma^{4}}{4 \theta^{2}} \bigg)
ds+\frac{\sigma^{4}}{4 \theta^{2}}\int_{1}^{t}(\widehat{\theta}_{s}-
\theta)^{2}ds.
\end{equation}
D'apr\`es la relation ({\normalfont{\ref{e39}}}), on a
\begin{eqnarray}\label{e15}
(\log t)^{-\frac{1}{2}}\int_{1}^{t}(\widehat{\theta}_{s}-
\theta)^{2}\bigg( \frac{ \langle M\rangle_{s}^{2}
}{s^{2}}-\frac{\sigma^{4}}{4 \theta^{2}} \bigg)
ds=~~~~~~~~~~~~~~~~~~~~~~~~~~~~~~~~~~~~~~~~~~\nonumber\\\nonumber\\(\log
t)^{-\frac{1}{2}} \int_{1}^{t}\mathcal{O}\bigg(\frac{\log\log
s}{s}\bigg)\bigg( \frac{ \langle M\rangle_{s}^{2}
}{s^{2}}-\frac{\sigma^{4}}{4 \theta^{2}} \bigg) ds \quad p.s..
\end{eqnarray}
L'hypoth\`ese (H2), implique que
$$ \frac{\langle M\rangle_{t}^{2}}{t^{2}}-\frac{\sigma^{4}}{4
\theta^{2}}= \mathbf{o}\bigg((\log\log t)^{-1}(\log
t)^{-\frac{1}{2}}\bigg)\quad p.s..$$ Sous cette derni\ère
hypoth\èse, on voit que
$$ (\log t )^{-\frac{1}{2}} \int_{1}^{t} (\widehat{\theta}_{s}-
\theta)^{2}\bigg( \frac{ \langle M\rangle_{s}^{2}
}{s^{2}}-\frac{\sigma^{4}}{4 \theta^{2}} \bigg)
ds\,\longrightarrow\,0,~~(t\longrightarrow\infty).$$
 Par suite
$$ \frac{\sigma^{2}}{2 \theta}(\log t
)^{-\frac{1}{2}} \int_{1}^{t}(\widehat{\theta}_{s}- \theta)^{2}ds-
\frac{\sigma^{2}}{2 \theta}(\log t )^{\frac{1}{2}} \Longrightarrow
\mathcal{N}(0,\frac{\sigma^{4}}{
 \theta^{2}}).$$
Ce qui signifie que  $$ (\log t )^{\frac{1}{2}}\bigg[ (2\log
t)^{-1}\int_{1}^{t}(\widehat{\theta}_{s}- \theta)^{2}ds-\theta
\bigg]\Longrightarrow\mathcal{N}(0,({2\theta})^{2}). $$
 Ainsi le r\'esultat est
\'etabli gr\^ace \à la convergence $(\ref{e34})$.\\
\item  En appliquant la loi du logarithme it\ér\é logarithmique au couple $(M,V)$ \\(voir th\éor\ème $5$ dans \cite{Fa}), on obtient $$ \limsup_{t\to \infty}(2\log t\log\log\log t)^{-\frac{1}{2}} \bigg|
 \int_{1}^{t}(\frac{M_{s}^{2}}{s}-\frac{\sigma^{2}}{2 \theta})\frac{ds}{s}\bigg|=\frac{\sigma^{2}}{\theta}\quad p.s..$$
 Gr\^ace \à la relation (\ref{e5}), on a $$\int_{1}^{t}\frac{M_{s}^{2}}{s}\,\frac{ds}{s}=\sigma^{-2}\int_{1}^{t}\frac{
\langle M\rangle_{s}^{2} }{s^{2}}(\widehat{\theta}_{s}- \theta)^{2}
ds= \sigma^{-2}I_{t}.$$ Vu les deux relations (\ref{e14}) et
(\ref{e15}), il vient que $$(2\log t\log\log\log
t)^{-\frac{1}{2}}\int_{1}^{t}(\widehat{\theta}_{s}-
\theta)^{2}\bigg( \frac{ \langle M\rangle_{s}^{2}
}{s^{2}}-\frac{\sigma^{4}}{4 \theta^{2}} \bigg)
ds\,\longrightarrow\,0,~~(t\longrightarrow\infty).$$ Ce qui implique
$$\limsup_{t\to \infty}(2\log t\log\log\log t)^{-\frac{1}{2}}\log t
\bigg| \frac{1}{2\log t}
\int_{1}^{t}\bigg(\widehat{\theta}_{s}-\theta)^{2}ds-\theta \bigg|=
2\theta \quad p.s..$$
 Ainsi compte tenu de la convergence (\ref{e34}), on d\éduit le r\ésultat.
\end{enumerate}
\subsection*{Preuve du th\'eor\`eme \ref{t4}}
D'apr\ès la loi du logarithme it\ér\é appliqu\ée \à la martingale
continue $\tilde{M}$ (voir \cite{PL}), on a
\begin{equation}{\label{e28}}
V_{t}^{-1}\tilde{M}_{t}=\mathcal{O}\bigg( (\log\log
V_{t})^{\frac{1}{2}}\bigg)\quad p.s..
\end{equation}
Vu la relation (\ref{e16}), la propri\ét\é ($P_{3})$ du lemme \ref{a} et la propri\ét\é ii) du\\
lemme \ref{b}, on conclut que
\begin{equation}{\label{e18}}
|\tilde{\theta}_{t}- \theta |= \mathcal{O}\bigg( (\frac{\log
t}{t^{\alpha}})^{\frac{1}{2}}\bigg)\quad p.s..
\end{equation}
D'o\`u la premi\ère assertion du th\éor\ème.\\
Par ailleurs gr\^ace aux relations (\ref{e16}) et (\ref{e27}), on
obtient
$$t^{\frac{1}{2}}(\bar{\theta}_{t}-\theta)=\sigma
t^{-\frac{1}{2}}\dint_{0}^{t}\dfrac{\tilde{M}_{s}}{P_{s}}\,ds.$$ La
propri\ét\é $(P_{4})$ implique que $$\dfrac{V_{s}}{U_{s}} =
s^{-\frac{\alpha}{2}}+o(s^{\frac{\alpha}{2}-1}))\quad p.s..$$
D'apr\ès ii) du lemme \ref{b}, on a
$$\dfrac{U_{s}}{P_{s}}=\frac{2\theta}{\sigma^{2}}+o(s^{\alpha'-\alpha})\quad p.s..$$
Par cons\équent vu que $\alpha'-\frac{3}{2}\alpha >
\alpha'-\frac{\alpha}{2}-1$ et $\alpha' > 2\alpha-1$, on obtient
$$t^{\frac{1}{2}}(\bar{\theta}_{t}-\theta)=\frac{2\theta}{\sigma}
t^{-\frac{1}{2}}\dint_{0}^{t}s^{-\frac{\alpha}{2}}\dfrac{\tilde{M}_{s}}{V_{s}}\,ds+
t^{-\frac{1}{2}}\dint_{0}^{t}o(s^{\alpha'-\frac{3}{2}\alpha})\dfrac{\tilde{M}_{s}}{V_{s}}\,ds.$$
Posons $Z_{s}=\dfrac{\tilde{M}_{s}}{V_{s}}$. D'apr\ès la relation
(\ref{e28}), on montre que
\begin{equation}{\label{e11}}
Z_{s}=\mathcal{O}\big((\log s)^{\frac{1}{2}}\big)\quad p.s..
\end{equation}
On d\éduit que
\begin{equation}{\label{e41}}
t^{\frac{1}{2}}(\bar{\theta}_{t}-\theta)=\frac{2\theta}{\sigma}
t^{-\frac{1}{2}}\dint_{0}^{t}s^{-\frac{\alpha}{2}}\dfrac{\tilde{M}_{s}}{V_{s}}\,ds+\mathcal{O}\bigg(t^{\alpha'-\frac{3}{2}\alpha+\frac{1}{2}}(\log
t)^{\frac{1}{2}}\bigg)\quad p.s..
\end{equation}
Comme  $\dfrac{dV_{s}}{V_{s}}\sim
\dfrac{s^{-\alpha}}{2}ds,~~(s\longrightarrow\infty).$ Alors
\begin{equation}{\label{e42}}
\frac{1}{2}\dint_{0}^{t}s^{-\frac{\alpha}{2}}Z_{s}ds=-\underbrace{\dint_{0}^{t}s^{\frac{\alpha}{2}}dZ_{s}}_{(K_{t})}+\underbrace{\dint_{0}^{t}\dfrac{s^{\frac{\alpha}{2}}}{V_{s}}d\tilde{M}_{s}}_{(L_{t})}.
\end{equation}
D'une part
$$K_{t}=t^{\frac{\alpha}{2}}Z_{t}-\frac{\alpha}{2}\dint_{0}^{t}s^{\frac{\alpha}{2}-1}Z_{s}ds.$$
Gr\^ace \à la relation (\ref{e11}), il vient que
\begin{equation}{\label{e45}}
K_{t}=\mathcal{O}\bigg((t^{\alpha}\log t )^{\frac{1}{2}}\bigg)\quad
p.s..
\end{equation}
 D'autre part, on a $$\langle
L\rangle_{t}=t^{\alpha}\dfrac{\langle
\tilde{M}\rangle_{t}}{V_{t}^{2}}-\alpha\dint_{0}^{t}s^{\alpha-1}\dfrac{\langle
\tilde{M}\rangle_{s}}{V_{s}^{2}}ds-2\dint_{0}^{t}s^{\alpha}\dfrac{\langle
\tilde{M}\rangle_{s}}{V_{s}^{2}}\dfrac{dV_{s}}{V_{s}}.$$ L'assertion
i) du lemme \ref{b} implique que $$\langle
L\rangle_{t}=\mathcal{O}(t)\quad p.s..$$ Encore une fois, la LLI
appliqu\ée \à la martingale continue $\tilde{M}$ montre que
\begin{equation}{\label{e43}}
L_{t}=\mathcal{O}\bigg((t\log\log t)^{\frac{1}{2}}\bigg)\quad p.s..
\end{equation}
En ins\érant (\ref{e45}) et (\ref{e43}) dans (\ref{e42}), on obtient
\begin{equation}{\label{e44}}
\dint_{0}^{t}s^{-\frac{\alpha}{2}}Z_{s}ds=\mathcal{O}\bigg((t\log\log
t)^{\frac{1}{2}}\bigg)+\mathcal{O}\bigg((t^{\alpha}\log t
)^{\frac{1}{2}}\bigg)\quad p.s..
\end{equation}
En combinant (\ref{e41}) et (\ref{e44}), la deuxi\ème assertion du
th\éo\èreme est \établie.
\subsection*{Preuve du th\'eor\`eme \ref{t5}}
\begin{enumerate}
\item
En appliquant le th\'eor\`eme de la limite centrale presque-s\^ure
pour le couple $(\tilde{M},V)$ avec
$V_{t}^{2}=\dint_{0}^{t}\omega_{s}^{2}ds$, on obtient
$$(\log V_{t}^{2})^{-1}\dint_{0}^{1}\delta_{\{V_{s}^{-1}\tilde{M}_{s}\}}d(\log V_{s}^{2})\Longrightarrow \mathcal{N}(0,\frac{\sigma^{2}}{2\theta})\quad p.s..$$
D'apr\ès la propri\ét\é ($P_{3})$ du lemme \ref{a}, on a
$$\frac{1-\alpha}{t^{1-\alpha}}\dint_{0}^{1}\delta_{\{V_{s}^{-1}\tilde{M}_{s}\}}\frac{ds}{s^{\alpha}}\Longrightarrow \mathcal{N}(0,\frac{\sigma^{2}}{2\theta})\quad p.s..$$
Gr\^ace \à la relation (\ref{e16}),  on obtient
$$\frac{1-\alpha}{t^{1-\alpha}}\dint_{0}^{t}\dfrac{ds}{s^{\alpha}}\,
\delta_{\{\frac{\sigma}{2\theta}s^{\frac{\alpha}{2}}(
\tilde{\theta}_{s}- \theta)\}}\Longrightarrow
\mathcal{N}(0,\frac{\sigma^{2}}{2\theta})\quad p.s..$$ Le r\ésultat
en d\écoule.\\
\item
i) \quad D'apr\ès la LFQ1 appliqu\ée \à la martingale $\tilde{M}$
normalis\ée par le processus
$V_{t}^{2}=\dint_{0}^{t}\omega_{s}^{2}ds$, on a
$$(\log
V_{t}^{2})^{-1}\dint_{0}^{t}\dfrac{\tilde{M_{s}^{2}}}{V_{s}^{2}}\dfrac{dV_{s}^{2}}{V_{s}^{2}}\longrightarrow\frac{\sigma^{2}}{2
\theta}\quad p.s.,~~(t\longrightarrow\infty).$$
 Vu la relation
(\ref{e16}), on a
\begin{equation}
\tilde{M}_{t}^{2}=\sigma^{-2}P_{t}^{2}({\tilde{\theta}_{t}-\theta})^{2}.
\end{equation}
Compte tenu de cette relation et de la propri\ét\é ($P_{3}$) du
lemme \ref{a}, on voit imm\édiatement que
\begin{equation}{\label{e17}}
\dfrac{1-\alpha}{ t^{1-\alpha}}\dint_{0}^{t}
\frac{P_{s}^{2}}{U_{s}^{2}}( \tilde{\theta}_{s}-
 \theta)^{2}ds \,
 \longrightarrow \dfrac{\sigma^{4}}{2 \theta}\quad p.s.,~~
 (t\longrightarrow\infty).
 \end{equation}
ii) \quad En appliquant la LFQ2 au couple $(\tilde{M},V)$, on
obtient
 $$(\log \langle \tilde{M}\rangle_{t})^{-1}\dint_{0}^{t}\frac{\tilde{M}_{s}^{2}}{ \langle \tilde{M} \rangle_{s}^{2} }d\langle \tilde{M} \rangle_{s}\longrightarrow 1\quad
 p.s.,~~(t\longrightarrow\infty).$$
Compte tenu de la propri\ét\é ($P_{3}$) du lemme \ref{a} et de la
propri\ét\é i) du\\ lemme \ref{b}, on a
$$\frac{1-\alpha}{t^{1-\alpha}}\dint_{0}^{t}\frac{\tilde{M}_{s}^{2}}{ \langle \tilde{M}
\rangle_{s}^{2} }d\langle \tilde{M} \rangle_{s}\longrightarrow
1\quad
 p.s.,~~(t\longrightarrow\infty).$$
Vu les relations (\ref{e4}) et (\ref{e16}), on obtient
$$\frac{\sigma^{2}(1-\alpha)}{t^{1-\alpha}}\dint_{0}^{t}\frac{{P}_{s}^{2}}{
\langle \tilde{M} \rangle_{s}^{2} }(\tilde{\theta}_{s}-
 \theta)^{2}\omega_{s}^{2}X_{s}^{2}ds\longrightarrow
1\quad
 p.s.,~~(t\longrightarrow\infty).$$
Gr\^ace \à la propri\ét\é $(P_{2}$) du lemme \ref{a}, on a
\begin{equation}{\label{e26}}
 \frac{4(1-\alpha)}{t^{1-\alpha}}\dint_{0}^{t} ( \tilde{\theta}_{s}-
 \theta)^{2}X_{s}^{2}ds \,\longrightarrow\, \sigma^{2}\quad
 p.s.,~~(t\longrightarrow\infty).
\end{equation}
Notons  que
\begin{eqnarray*}
   \bigg|\dfrac{4(1-\alpha)}{t^{1-\alpha}}\dint_{0}^{t}(\tilde{\theta}_{s}- \theta
)^{2}X^{2}_{s}ds-\frac{4(1-\alpha)}{t^{1-\alpha}}\dint_{0}^{t}(\tilde{\theta}_{s}-
\bar{\theta}_{t})^{2}X^{2}_{s}ds\bigg|\leq~~~~~~~~~~~~~~~~~~~~~~~~\\\\
\underbrace{\frac{8(1-\alpha)}{t^{1-\alpha}}(\bar{\theta}_{t}-
\theta )\int_{0}^{t}(\tilde{\theta}_{s}- \bar{\theta}_{t}
)X^{2}_{s}ds}_{(E_{t})}
+\underbrace{\frac{4(1-\alpha)}{t^{1-\alpha}}(\bar{\theta}_{t}-
\theta )^{2}\int_{0}^{t}X^{2}_{s}ds}_{(F_{t})}.
\end{eqnarray*}
D'apr\ès le th\éor\ème \ref{t4} et la relation ({\ref{e9}}), on
obtient
\begin{equation}{\label{e29}}
E_{t}\longrightarrow 0 \quad p.s.\quad \mbox{et}\quad
F_{t}\longrightarrow 0 \quad p.s.,~~(t\longrightarrow\infty).
\end{equation}
 Gr\^ace \à la
convergence (\ref{e26}), le r\ésultat est \établi.\\\\
 iii)\quad Posons $$\tilde{I_{t}}=\dfrac{1-\alpha}{
t^{1-\alpha}}\dint_{0}^{t} \frac{P_{s}^{2}}{U_{s}^{2}}(
\tilde{\theta}_{s}-
 \theta)^{2}ds.$$
Soit
$$\tilde{I_{t}}=\tilde{I_{t}^{1}}+\tilde{I_{t}^{2}}$$
avec $$\begin{array}{ccl}
      \tilde{I_{t}^{1}} & = &  \dfrac{1-\alpha}{
t^{1-\alpha}}\dint_{0}^{t} ( \tilde{\theta}_{s}-
 \theta)^{2}(\frac{P_{s}^{2}}{U_{s}^{2}}-\dfrac{\sigma^{4}}{4
 \theta^{2}})ds\quad\mbox{et}\\\\
      \tilde{I_{t}^{2}} & = & \dfrac{\sigma^{4}(1-\alpha)}{4 \theta^{2}t^{1-\alpha}}\dint_{0}^{t}(\tilde{\theta}_{s}-
 \theta)^{2} ds.
    \end{array}$$
 Vu l'hypoth\èse $(H3)$, la relation  ii) du lemme \ref{b}
 implique que
\begin{equation}\label{e19}
\bigg(\int_{0}^{t}\omega_{s}ds \bigg)^{-2}
\bigg(\int_{0}^{t}\omega_{s}X_{s}^{2}ds
\bigg)^{2}-\frac{\sigma^{4}}{4\theta^{2}}=
\mathbf{o}\bigg(t^{2(\alpha-1)}(\log t)^{-1}\bigg)\quad p.s..
\end{equation}
D'apr\ès le th\éor\ème \ref{t4} et la relation (\ref{e19}), il vient
que
$$\bigg|\tilde{I_{t}}-\tilde{I_{t}^{2}}\bigg|\longrightarrow 0, \quad
(t\longrightarrow\infty).$$ Gr\^ace \à la propri\ét\é (\ref{e17}),
on obtient
\begin{equation}\label{e36} \frac{1-\alpha}{
2t^{1-\alpha}}\int_{0}^{t} (\tilde{\theta}_{s}-\theta)^{2}ds
\longrightarrow  \theta\quad
 ~p.s.,~~(t\longrightarrow\infty).
\end{equation}
Par ailleurs, pour $\dfrac{1}{2}<\alpha<1$, on voit que\\
\begin{eqnarray*}
\Bigg|\dfrac{1-\alpha}{ 2t^{1-\alpha}}\dint_{0}^{t}
(\tilde{\theta}_{s}-\theta)^{2}ds-\dfrac{1-\alpha}{
2t^{1-\alpha}}\dint_{0}^{t}
(\tilde{\theta}_{s}-\bar{\theta}_{t})^{2}ds\Bigg|\leq~~~~~~~~~~~~~~~~~~~~~~~~~~~~~~~~~~~~~\\\\
\underbrace{\dfrac{1-\alpha}{
t^{1-\alpha}}(\bar{\theta}_{t}-\theta)\int_{0}^{t}
(\tilde{\theta}_{s}-\bar{\theta}_{t})ds}_{(Q_{t})}+\underbrace{\dfrac{1-\alpha}{
2t^{-\alpha}}(\bar{\theta}_{t}-\theta)^{2}}_{(S_{t})}.
\end{eqnarray*}
\\Le th\éor\ème \ref{t4} implique que
$$Q_{t}\longrightarrow0\quad
 p.s.\quad \mbox{et}\quad S_{t}\longrightarrow0\quad
 p.s.,~~(t\longrightarrow\infty).$$
Par cons\équent
\begin{equation}\label{e37}
\Bigg|\dfrac{1-\alpha}{ 2t^{1-\alpha}}\dint_{0}^{t}
(\tilde{\theta}_{s}-\theta)^{2}ds-\dfrac{1-\alpha}{
2t^{1-\alpha}}\dint_{0}^{t}
(\tilde{\theta}_{s}-\bar{\theta}_{t})^{2}ds\Bigg|\longrightarrow
0\quad
 ~p.s.,~~(t\longrightarrow\infty).
\end{equation}
Alors compte tenu de la convergence $(\ref{e36})$, le r\ésultat est
\établi.
\end{enumerate}
\subsection*{Preuve du th\'eor\`eme \ref{t6}}
\begin{enumerate}
\item
D'apr\ès i) du lemme \ref{b} et vu que $\dfrac{1}{2}\leq \alpha'<
3\alpha-2$, l'hypoth\èse $(H4)$ implique que
$$ \bigg(\int_{0}^{t}\omega_{s}^{2}ds
\bigg)^{-1} \bigg(\int_{0}^{t}\omega_{s}^{2}X_{s}^{2}ds
\bigg)-\frac{\sigma^{2}}{2\theta}= \mathbf{o}(t^{2(\alpha-1)})\quad
p.s..$$ Alors en appliquant le TLCL au couple
$(\tilde{M},V)$ pour $V_{t}^{2}=\dint_{0}^{t}\omega_{s}^{2}\,ds $\\
et $f(x)=x^{2}-1$, on obtient
\begin{equation}{\label{e22}}
(\log
V_{t}^{2})^{-\frac{1}{2}}\dint_{1}^{t}(\frac{\tilde{M_{s}^{2}}}{V_{s}^{2}}-\frac{\sigma^{2}}{2
\theta})\frac{dV_{s}^{2}}{V_{s}^{2}}\Longrightarrow\mathcal{N}(0,\frac{\sigma^{4}}{\theta^{2}}).
\end{equation}
Gr\^ace \à la propri\ét\é ($P_{3}$) du lemme \ref{a} et la relation
(\ref{e16}), on obtient
\begin{equation}{\label{e21}}
(\frac{t^{1-\alpha}}{1-\alpha})^{-\frac{1}{2}}\dint_{1}^{t}\frac{\tilde{M}_{s}}{V_{s}^{2}}\frac{ds}{s^{\alpha}}\sim
\sigma^{-2}(\frac{t^{1-\alpha}}{1-\alpha})^{-\frac{1}{2}}\dint_{1}^{t}\dfrac{P_{s}^{2}}{U_{s}^{2}}({\tilde{\theta}_{t}-\theta})^{2}ds,~~(t\longrightarrow\infty).
\end{equation}
Soit $\tilde{J_{t}}$ le terme de droite de la derni\ère
\équivalence. On pose
$$\tilde{J_{t}}=\tilde{J_{t}^{1}}+\tilde{J_{t}^{2}}$$
avec $$\begin{array}{rcl}
  \tilde{J_{t}^{1}}&=&\sigma^{-2}(\dfrac{t^{1-\alpha}}{1-\alpha})^{-\frac{1}{2}}\dint_{1}^{t}({\tilde{\theta}_{t}-\theta})^{2}(\dfrac{P_{s}^{2}}{U_{s}^{2}}-\frac{
\sigma^{4} }{4 \theta^{2}})ds \quad \mbox{et} \\\\
  \tilde{J_{t}^{2}} & = & \dfrac{
\sigma^{2}}{4
\theta^{2}}(\dfrac{t^{1-\alpha}}{1-\alpha})^{-\frac{1}{2}}\dint_{1}^{t}({\tilde{\theta}_{t}-\theta})^{2}ds.
\end{array}$$
Vu l'hypoth\èse $(H4)$ et d'apr\ès ii) du lemme \ref{b}, il vient
que pour $\; \alpha'< 3\alpha-2$
\begin{equation}\label{e20}
\bigg(\int_{0}^{t}\omega_{s}ds \bigg)^{-2}
\bigg(\int_{0}^{t}\omega_{s}X_{s}^{2}ds
\bigg)^{2}-\frac{\sigma^{4}}{4\theta^{2}}=
\mathbf{o}\bigg(t^{\frac{(\alpha-1)}{2}}(\log t)^{-1}\bigg)\quad
p.s..
\end{equation}
D'apr\ès le th\éor\ème \ref{t4} et la relation (\ref{e20}), on
d\éduit que
\begin{equation}\label{e23}
\bigg|\tilde{J_{t}}-\tilde{J_{t}^{2}}\bigg|\longrightarrow 0, \quad
(t\longrightarrow\infty).
\end{equation}
 En combinant  (\ref{e22}), (\ref{e21}) et (\ref{e23}), on obtient
$$\frac{\sigma^{2}}{4\theta^{2}}(\frac{t^{1-\alpha}}{1-\alpha})^{-\frac{1}{2}}\dint_{0}^{t}(\tilde{\theta}_{s}-
\theta)^{2}ds-\frac{\sigma^{2}}{2\theta}(\frac{t^{1-\alpha}}{1-\alpha})^{\frac{1}{2}}\Longrightarrow
\mathcal{N}(0,\frac{\sigma^{4}}{\theta^{2}}).$$ Alors
$$(\frac{t^{1-\alpha}}{1-\alpha})^{\frac{1}{2}}\Bigg
(\dfrac{1-\alpha}{ 2t^{1-\alpha}}\dint_{0}^{t} (\tilde{\theta}_{s}-
\theta)^{2}ds-\theta \Bigg) \Longrightarrow
\mathcal{N}(0,4\theta^{2}).$$ Gr\^ace
\à la convergence $(\ref{e37})$, on obtient le r\ésultat.\\
\item
En appliquant la LLIL au couple $(\tilde{M},V)$, on obtient $$
\limsup_{t\to \infty}(2\log V_{t}^{2}\log\log\log
V_{t}^{2})^{-\frac{1}{2}}
\bigg|\dint_{1}^{t}(\frac{\tilde{M_{s}^{2}}}{V_{s}^{2}}-\frac{\sigma^{2}}{2
\theta})\frac{dV_{s}^{2}}{V_{s}^{2}}
 \bigg|=\frac{\sigma^{2}}{\theta}\quad p.s..$$
 Ainsi d'apr\ès (\ref{e23}), on a $$\dint_{1}^{t}\frac{\tilde{M_{s}^{2}}}{V_{s}^{2}}\frac{ds}{{s}^{\alpha}}\sim(\frac{t^{1-\alpha}}{1-\alpha})^{\frac{1}{2}}\tilde{J_{t}^{2}},~~(t\longrightarrow\infty).$$
Gr\^ace \à cette derni\ère \équivalence et la propri\ét\é $(P_{3}$)
du lemme \ref{a}, on obtient
$$\limsup_{t\to \infty}(2\log V_{t}^{2}\log\log\log
V_{t}^{2})^{-\frac{1}{2}}\dfrac{ t^{1-\alpha}}{1-\alpha}\Bigg
|\dfrac{1-\alpha}{ 2t^{1-\alpha}}\dint_{0}^{t} (\tilde{\theta}_{s}-
\theta)^{2}ds- \theta \Bigg|= 2\theta \quad p.s.. $$ Ce qui signifie
que 
$$\dlimsup_{t\to \infty} {\dfrac{
t^{1-\alpha}}{{\sqrt{\log\log t^{1-\alpha}}}}}\,\Bigg
|\dfrac{1-\alpha}{ 2t^{1-\alpha}}\dint_{0}^{t} (\tilde{\theta}_{s}-
\theta)^{2}ds- \theta \Bigg|= 2\theta \sqrt{2(1-\alpha)}\quad
p.s..$$ Le r\ésultat est \établi vu la convergence $(\ref{e37})$. Ce
qui ach\ève la preuve du th\éor\ème.
\end{enumerate}
\section{Annexe}
\subsection*{Preuve du lemme \ref{a}}
\begin{enumerate}
\item {\bf Preuve de ($P_{1}$)}\\
Soit $N=(N_{t},\,t\geq0)$ la fonction d\éfinie par
$$N_{t}=t^{-\alpha}\omega_{t}^{-1}U_{t}.$$
D'apr\ès l'expression du poids $(\omega_{t})$, on voit que
$$N_{t}=t^{-\frac{\alpha}{2}}e^{-\frac{t^{1-\alpha}}{2(1-\alpha)}}\dint_{0}^{t}s^{-\frac{\alpha}{2}}e^{\frac{s^{1-\alpha}}{2(1-\alpha)}}ds.$$
Dans la suite on s'interessera au comportement asymptotique de
 $(N_{t},\,t\geq 0)$. Il est clair que
  $$N_{t}\geq
t^{-\frac{\alpha}{2}}e^{-\frac{t^{1-\alpha}}{2(1-\alpha)}}\bigg(\dint_{t'}^{t}s^{-\frac{\alpha}{2}}e^{\frac{s^{1-\alpha}}{2(1-\alpha)}}ds\bigg),~~\mbox{avec}\quad
1\leq t'<t.$$ Donc $$N_{t}\geq
e^{-\frac{t^{1-\alpha}}{2(1-\alpha)}}\dint_{t'}^{t}s^{-{\alpha}}e^{\frac{s^{1-\alpha}}{2(1-\alpha)}}ds.$$
Comme
$$\dint_{t'}^{t}s^{-{\alpha}}e^{\frac{s^{1-\alpha}}{2(1-\alpha)}}ds=2\bigg(e^{\frac{t^{1-\alpha}}{2(1-\alpha)}}-e^{\frac{(t')^{1-\alpha}}{2(1-\alpha)}}\bigg).$$
D'o\`u, on a $$N_{t}\geq
2-2e^{\frac{(t')^{1-\alpha}-t^{1-\alpha}}{2(1-\alpha)}}.$$ Pour
$t'<t$, on d\éduit que $$N_{t}-2\geq
-2e^{\frac{(t')^{1-\alpha}-t^{1-\alpha}}{2(1-\alpha)}}=\mathbf{o}(t^{\alpha-1}).$$
Par cons\équent
 $$\liminf_{t\to
\infty}t^{1-\alpha}(N_{t}-2)=0.$$ D'o\`u
\begin{equation}\label{l3}
\liminf_{t\to \infty}N_{t}=2.
\end{equation}
Par ailleurs, notons que pour $0<\lambda<1$, on a
\begin{equation}\label{l4}
N_{t}\leq\beta_{t}^{\lambda}+\gamma_{t}^{\lambda}
\end{equation}
o\`u$$\begin{array}{ccc}
                      \beta_{t}^{\lambda} & = &t^{-\frac{\alpha}{2}}
                      e^{-\frac{t^{1-\alpha}}{2(1-\alpha)}}\dint_{1}^{t\lambda}s^{-{\alpha}}e^{\frac{s^{1-\alpha}}{2(1-\alpha)}}ds,
                      \\\\
                      \gamma_{t}^{\lambda} & = &t^{-\frac{\alpha}{2}}(t\lambda)^{\frac{\alpha}{2}}
                      e^{-\frac{t^{1-\alpha}}{2(1-\alpha)}}\dint_{t\lambda}^{t}s^{-{\alpha}}e^{\frac{s^{1-\alpha}}{2(1-\alpha)}}ds.
                    \end{array}$$
Afin de trouver la limite sup\érieure de $N_{t}$, on cherche le
comportement asymptotique des fonctions $\beta_{t}^{\lambda}$ et
$\gamma_{t}^{\lambda}$. En remarquant que $$\begin{array}{ccl}
   t^{1-\alpha}\beta_{t}^{\lambda}& = & 2t^{1-\frac{3}{2}\alpha}e^{-\frac{t^{1-\alpha}}{2(1-\alpha)}}\bigg(e^{-\frac{(t\lambda)^{1-\alpha}}{2(1-\alpha)}}-e^{-\frac{1}{2(1-\alpha)}}\bigg)
   \\\\
   & =
   &2t^{1-\frac{3}{2}\alpha}e^{-\frac{1}{2(1-\alpha)}\bigg((t\lambda)^{1-\alpha}-t^{1-\alpha}\bigg)}-2t^{1-\frac{3}{2}\alpha}e^{-\frac{1}{2(1-\alpha)}(1-t^{1-\alpha})}
\end{array}$$
et en utilisant le fait que $\frac{1}{2}<\alpha<1$ et $0<\lambda<1$,
on d\éduit que
\begin{equation}\label{l5}
t^{1-\alpha}\beta_{t}^{\lambda}\longrightarrow 0,
\quad(t\longrightarrow\infty).
\end{equation}
De la m\^eme fa\c con on \établit la relation suivante
$$\gamma_{t}^{\lambda}=2\lambda^{\frac{\alpha}{2}}\bigg(1-e^{-\frac{t^{1-\alpha}}{2(1-\alpha)}(\lambda^{1-\alpha}-1)}\bigg).$$
Comme $\lambda^{1-\alpha}-1<0$, on obtient$$\gamma_{t}^{\lambda}\sim
2\lambda^{\frac{\alpha}{2}},~~(t\longrightarrow\infty).$$ De plus on
a
\begin{equation}\label{l6}
t^{1-\alpha}\big(\gamma_{t}^{\lambda}-2\lambda^{\frac{\alpha}{2}}\big)=-t^{1-\alpha}e^{\frac{t^{1-\alpha}}{2(1-\alpha)}(\lambda^{1-\alpha}-1)}=\mathbf{o}(1).
\end{equation}
En rempla\c cant $\lambda$ par $\lambda_{t}=1-t^{-\alpha}$, les
propri\ét\és $(\ref{l4})$, $(\ref{l5})$ et $(\ref{l6})$ restent
vraies. Par cons\équent on obtient
\begin{equation}\label{l7}
t^{1-\alpha}(N_{t}-2)\leq
t^{1-\alpha}\beta_{t}^{\lambda_{t}}+t^{1-\alpha}\bigg(\gamma_{t}^{\lambda_{t}}-2\lambda_{t}^{\frac{\alpha}{2}}\bigg)+2t^{1-\alpha}\bigg(\lambda_{t}^{\frac{\alpha}{2}}-1\bigg).
\end{equation}
 Vu que $\lambda_{t}$ v\érifie
$$t^{1-\alpha}\big(1-\lambda_{t}\big)\longrightarrow 0, \quad(t\longrightarrow\infty)$$
on d\éduit que $$0\leq \lambda_{t}^{\frac{\alpha}{2}}-1\leq
\frac{\alpha}{2}(\lambda_{t}-1)\lambda_{t}^{\frac{\alpha}{2}-1}=\mathbf{o}(t^{\alpha-1}).$$
Compte tenu de ce dernier r\ésultat et le fait que\\
$$t^{1-\alpha}\beta_{t}^{\lambda_{t}}\longrightarrow 0 \quad
\mbox{et}\quad
t^{1-\alpha}\big(\gamma_{t}^{\lambda_{t}}-2\lambda_{t}^{\frac{\alpha}{2}}\big)\longrightarrow
0, \quad(t\longrightarrow\infty),$$
 il vient que
\begin{equation}\label{l8}
\limsup_{t\to \infty}N_{t}=2.
\end{equation}
Gr\^ace \à (\ref{l3}) et (\ref{l8}),  on conclut le r\ésultat.\\
\item {\bf Preuve de\,$(P_{2})$}\\
Posons $N'=(N'_{t},\,t\geq0)$, la fonction d\éfinie par
$$N'_{t}=t^{-\alpha}\omega_{t}^{-2}V_{t}^{2}.$$
D'apr\ès l'expression du poids $(\omega_{t})$, on voit que
$$N'_{t}=e^{-\frac{t^{1-\alpha}}{1-\alpha}}\dint_{0}^{t}s^{-\alpha}e^{\frac{s^{1-\alpha}}{1-\alpha}}ds=1-e^{-\frac{t^{1-\alpha}}{1-\alpha}}\longrightarrow 1,
\quad(t\longrightarrow\infty).$$
 Par ailleurs
$$t^{1-\alpha}(N'_{t}-1)=-t^{1-\alpha}e^{-\frac{t^{1-\alpha}}{1-\alpha}}\longrightarrow 0,
\quad(t\longrightarrow\infty).$$ D'o\`u  la propri\ét\é $(P_{2})$.\\
\item {\bf Preuve de\,$(P_{3})$}\\
La propri\ét\é $(P_{2})$ implique que
$$t^{-\alpha}=\dfrac{\omega_{t}^{2}}{V_{t}^{2}}+\mathbf{o}(\dfrac{\omega_{t}^{2}}{V_{t}^{2}}t^{\alpha-1}).$$
Compte tenu de  l'expression du poids, on obtient
$$V_{t}^{2}=e^{\frac{t^{1-\alpha}}{1-\alpha}}-1 \sim t^{\alpha}\omega_{t}^{2}, \quad(t\longrightarrow\infty).$$
 Par cons\équent on a
$$t^{-\alpha}-V_{t}^{-2}\omega_{t}^{2}=\mathbf{o}(t^{-1}).$$
Vu que
$$\dfrac{\omega_{t}^{2}}{V_{t}^{2}}\sim V_{t}^{-2}\dfrac{dV_{t}^{2}}{dt}, \quad(t\longrightarrow\infty).$$
La propri\ét\é $(P_{3})$ est \établie.\\
\item {\bf Preuve de\,$(P_{4})$}\\
Gr\^ace \à la propri\ét\é $(P_{1})$, on a
$$t^{-2\alpha}\omega_{t}^{-2}U_{t}^{2}=4+\mathbf{o}(t^{\alpha-1}).$$
La propri\ét\é $(P_{2})$ donne
$$t^{-\alpha}\omega_{t}^{-2}V_{t}^{2}=\frac{1}{4}+\mathbf{o}(t^{\alpha-1}).$$
En combinant ces deux r\ésultats, on obtient
$$\frac{V_{t}^{2}}{U_{t}^{2}}=t^{-\alpha}+\mathbf{o}(t^{-1}).$$
Ce qui ach\ève la preuve de $(P_{4})$.
\end{enumerate}
\subsection*{Preuve du lemme \ref{b}}
D'apr\ès l'\égalit\é (\ref{e3}), on a$$\langle
\tilde{M}\rangle_{t}=\frac{\sigma^{2}}{2\theta}V_{t}^{2}+\int_{0}^{t}\omega_{s}^{2}(X_{s}^{2}-\frac{\sigma^{2}}{2\theta})ds.$$
Comme$$\frac{\langle
{M}\rangle_{t}}{t}=\frac{1}{t}\int_{0}^{t}X_{s}^{2}ds\longrightarrow
\frac{\sigma^{2}}{2\theta}~p.s.,~~ (t\longrightarrow\infty).$$ Alors
$$A_{t}-\frac{\sigma^{2}}{2\theta}=\frac{1}{t}\int_{0}^{t}(X_{s}^{2}-\frac{\sigma^{2}}{2\theta})ds,$$
 avec $A_{t}=\frac{\langle {M}\rangle_{t}}{t}~,~t>0$.
Donc
$$d(t(A_{t}-\frac{\sigma^{2}}{2\theta}))=(X_{t}^{2}-\frac{\sigma^{2}}{2\theta})dt.$$
Par suite
$$\begin{array}{rcl}
  \dint_{0}^{t}\omega_{s}^{2}(X_{s}^{2}-\frac{\sigma^{2}}{2\theta})ds  & = & \dint_{0}^{t}\omega_{s}^{2}d(s(A_{s}-\frac{\sigma^{2}}{2\theta})) \\
   & =
   &[\omega_{s}^{2}s(A_{s}-\frac{\sigma^{2}}{2\theta})]_{0}^{t}-\dint_{0}^{t}s(A_{s}-\frac{\sigma^{2}}{2\theta})d\omega_{s}^{2}.
\end{array}$$
Compte tenu de l'expression de la pond\ération d\éfinie dans
(\ref{e8}), on voit que
$$d\omega_{s}^{2}=\omega_{s}^{2}\{\frac{1}{(1-\alpha)s^{\alpha}}-\frac{\alpha}{s}\}ds.$$
Par cons\équent
\begin{equation}\label{e38}
 \langle
\tilde{M}\rangle_{t}-\frac{\sigma^{2}}{2\theta}V_{t}^{2}=[\omega_{s}^{2}s(A_{s}-\frac{\sigma^{2}}{2\theta})]_{0}^{t}+\alpha\dint_{0}^{t}(A_{s}-\frac{\sigma^{2}}{2\theta})\omega_{s}^{2}ds-\frac{1}{1-\alpha}\dint_{0}^{t}s^{1-\alpha}(A_{s}-\frac{\sigma^{2}}{2\theta})\omega_{s}^{2}ds.
\end{equation} Gr\^ace \à la propri\ét\é $(P_{2}$) et l'hypoth\èse $(H3)$, il vient que
  $$\dfrac{1}{V_{t}^{2}}\dint_{0}^{t}(A_{s}-\frac{\sigma^{2}}{2\theta})\omega_{s}^{2}ds  \leq  \dint_{0}^{t}(A_{s}-\frac{\sigma^{2}}{2\theta})\dfrac{\omega_{s}^{2}}{V_{s}^{2}}ds
  =\mathbf{o}(t^{\alpha'-\alpha}).$$
De m\^eme on a
$$\dfrac{1}{V_{t}^{2}}\dint_{0}^{t}s^{1-\alpha}(A_{s}-\frac{\sigma^{2}}{2\theta})\omega_{s}^{2}ds=\mathbf{o}(t^{\alpha'-\alpha}).$$
Vu la relation $(\ref{e38})$, on obtient la propri\ét\é i) suivante
$$\dfrac{\langle
\tilde{M}\rangle_{t}}{V_{t}^{2}}-\frac{\sigma^{2}}{2\theta}=\mathbf{o}(t^{\alpha'-\alpha})~~p.s..$$
De la m\^eme fa\c con, on \établit la propri\ét\é ii).
 \bibliographystyle{unsrt}

\end{document}